\documentclass{amsart}

\usepackage{fancyhdr}
\usepackage{lastpage}
\usepackage{stmaryrd}

\newcommand{\bdis}{\begin{displaymath}}
\newcommand{\edis}{\end{displaymath}}
\newcommand{\be}{\begin{equation}}
\newcommand{\ee}{\end{equation}}

\newcommand{\mcal}{\mathcal}

\theoremstyle{definition}

\theoremstyle{remark}
\newtheorem{remark}[]{Remark}

\newtheorem*{mydef1}{{\bf Theorem}}

\newtheorem*{mydef4}{{\bf Corollary}}

\newtheorem*{mydef5}{{\bf Lemma}}

\numberwithin{equation}{section}

\begin{document}

\title{Jacob's ladders and the first asymptotic formula for the expression of the sixth order
$|\zeta(1/2+i\varphi(t)/2)|^4|\zeta(1/2+it)|^2$}

\author{Jan Moser}

\address{Department of Mathematical Analysis and Numerical Mathematics, Comenius University, Mlynska Dolina M105, 842 48 Bratislava, SLOVAKIA}

\email{jan.mozer@fmph.uniba.sk}


\begin{abstract}
It is proved in this paper that there is a fine correlation between the values of $|\zeta(1/2+i\varphi(t)/2)|^4$ and
$|\zeta(1/2+it)|^2$ which correspond to two segments with gigantic distance each from other. This new asymptotic formula
cannot be obtained in known theories of Balasubramanian, Heath-Brown and Ivic. \\

Dedicated to the memory of Anatolij Alekseevich Karatsuba (1937-2008).
\end{abstract}

\maketitle

\section{Results}

\subsection{}

The following theorem holds true

\begin{mydef1}
\be \label{f1.1}
\int_T^{T+U}Z^4\left[\frac{\varphi(t)}{2}\right]Z^2(t){\rm d}t\sim \frac{1}{2\pi^2}U\ln^5 T,\ U=T^{7/8+2\epsilon},\ T\to\infty ,
\ee
where $\varphi(t)$ is a Jacob's ladder,
\be \label{f1.2}
t-\frac{1}{2}\varphi(t)\sim (1-c)\pi(t),\ t\to \infty
\ee
(see \cite{3}, (6.2)), $c$ is the Euler constant and $\pi(t)$ is the prime-counting function.
\end{mydef1}

\begin{remark}

The formula (\ref{f1.1}) is the first asymptotic formula in the theory of the Riemann zeta-function for the sixth order expression
$|\zeta(1/2+i\varphi(t)/2)|^4|\zeta(1/2+it)|^2$. This formula cannot be obtained by methods of Balasubramanian, Heath-Brown and Ivic (see \cite{2}).

\end{remark}

\subsection{}

Since (see (\ref{f1.2})
\bdis
T+U-\frac{1}{2}\varphi(T+U)>(1-c-\epsilon)\frac{T}{\ln T} ,
\edis
we have (see the condition for $U$ in (\ref{f1.1})
\be \label{f1.3}
T-\frac{1}{2}\varphi(T+U)>(1-c-\epsilon)\frac{T}{\ln T}-U>(1-c-2\epsilon)\frac{T}{\ln T} .
\ee
Hence, using the mean-value theorem in (\ref{f1.1}) we obtain

\begin{mydef4}

\be \label{f1.4}
Z^4\left[\frac{\varphi(\alpha)}{2}\right]Z^2(\alpha)\sim \frac{1}{2\pi^2}\ln^5 T,\ T\to \infty ,
\ee
where
\be \label{f1.5}
\alpha\in (T,T+U),\ \frac{1}{2}\varphi(\alpha)\in \left(\frac{1}{2}\varphi(T),\frac{1}{2}\varphi(T+U)\right), \
\alpha=\alpha(T,U,\varphi)=\alpha(T),
\ee
and
\begin{eqnarray} \label{f1.6}
& &
[T,T+U]\bigcap \left[\frac{1}{2}\varphi(T),\frac{1}{2}\varphi(T+U)\right]=\emptyset , \nonumber \\
& &
\rho\left\{ [T,T+U]; \left[\frac{1}{2}\varphi(T),\frac{1}{2}\varphi(T+U)\right]\right\}>(1-c-2\epsilon)\frac{T}{\ln T}\to \infty ,
\end{eqnarray}
where $\rho$ denotes the distance of the corresponding segments.

\end{mydef4}

\begin{remark}

Since $U=T^{7/8+\epsilon}$, $\alpha(T)$ is the single-valued function of $T$: one mean-value of the set $\{\alpha(T)\}$
corresponds to each sufficiently big $T$.

\end{remark}

\begin{remark}

Some \emph{nonlocal interaction} of the functions
\bdis
|\zeta(1/2+i\varphi(t)/2)|^4,\quad |\zeta(1/2+it)|^2
\edis
is expressed by formula (\ref{f1.4}). This interaction is connected with two segments unboundedly receding each from other (see (\ref{f1.6}),
$\rho\to\infty$ as $T\to\infty$) - like mutually receding galaxies (the Hubble law).

\end{remark}

\begin{remark}

Since $T\sim \alpha(T),\ \alpha\in (T,T+U)$, then from (\ref{f1.4}) (see also (\ref{f1.2}) we obtain
\begin{eqnarray} \label{f1.7}
& &
|Z[\alpha(T)]|\sim \frac{1}{\sqrt{2}\pi}\frac{\ln^{5/2}\alpha}{Z^2\left[\frac{\varphi(\alpha)}{2}\right]}, \nonumber \\
& &
\alpha(T)\sim\frac{1}{2}\varphi(\alpha)+(1-c)\pi(T) ,\ T\to\infty .
\end{eqnarray}
By (\ref{f1.7}) we have the prediction of the value $|Z(\alpha)|,\ \alpha\in (T,T+U)$  by means of the value $Z^2[\varphi(\alpha)/2]$ which
corresponds to argument $\varphi(\alpha)/2$ running in the segment $(1/2\varphi(T),1/2\varphi(T+U))$ which descended from very deep past
(see (\ref{f1.6})).

\end{remark}

\begin{remark}
The following auto-correlation formula
\be \label{f1.8}
\int_T^{T+U}Z^2\left[\frac{\varphi(t)}{2}\right]Z^2(t){\rm d}t\sim U\ln^2 T,\ U=T^{1/3+2\epsilon}
\ee
of the fourth order takes place. This formula is fully similar to correlation formula (\ref{f1.1}).

\end{remark}

\begin{remark}

The notions \emph{correlation, auto-correlation} and \emph{prediction} for the signal
\bdis
Z(t)=e^{i\theta(t)}\zeta\left(\frac{1}{2}+it\right)
\edis
generated by the Riemann zeta-function are inspired by the fundamental book of N. Wiener, \cite{8}. Incidentally, this book - for use of communication
engineering - contains multitude of stimuli for development of the theory of the zeta-function.
\end{remark}

This paper is a continuation of the series of papers \cite{3}-\cite{6}.

\section{$\hat{Z}^2$-transformation}

We start with the formula (see \cite{3}, (3.5), (3.9))
\be \label{f2.1}
Z^2(t)=\Phi^\prime_\varphi[\varphi(t)]\frac{{\rm d}\varphi}{{\rm d}t},\ T\geq T_0[\varphi],
\ee
where (see \cite{5}, (1.5))
\be \label{f2.2}
\Phi^\prime_\varphi[\varphi(t)]=\frac{1}{2}\left\{ 1+\mcal{O}\left( \frac{\ln \ln t}{\ln t}\right) \right\}\ln t .
\ee
Next we have
\be \label{f2.3}
\hat{Z}^2(t)=\frac{{\rm d}\varphi(t)}{{\rm d}t},\ t\in [T,T+U],\ T\geq T_0[\varphi],\ U\in \left(\left. 0,\frac{T}{\ln T}\right]\right.
\ee
by (\ref{f2.1}) (\ref{f2.2}), where
\be \label{f2.4}
\hat{Z}^2(t)=\frac{Z^2(t)}{\Phi^\prime_\varphi[\varphi(t)]}=\frac{2Z^2(t)}{\left\{ 1+\mcal{O}\left( \frac{\ln \ln t}{\ln t}\right) \right\}\ln t} .
\ee
Then we obtain from (\ref{f2.3}) the following lemma.

\begin{mydef5}

For every integrable function $f(x),\ x\in [1/2\varphi(T),1/2\varphi(T+U)]$ the following is true
\be \label{f2.5}
\int_T^{T+U}f\left[\frac{\varphi(t)}{2}\right]\hat{Z}^2(t){\rm d}t=2\int_{\frac{1}{2}\varphi(T)}^{\frac{1}{2}\varphi(T+U)}f(x){\rm d}x,\
U\in \left.\left( 0,\frac{T}{\ln T}\right.\right] .
\ee
\end{mydef5}

For example, in the case $f(x)=|Z(x)|^\Delta$ with $\Delta\in (0,\infty)$ we have the following  $\hat{Z}^2$-transformation
\be \label{f2.6}
\int_T^{T+U}\left| Z\left[\frac{\varphi(t)}{2}\right]\right|^\Delta \hat{Z}^2(t){\rm d}t=2\int_{\frac{1}{2}\varphi(T)}^{\frac{1}{2}\varphi(T+U)}
|Z(x)|^\Delta{\rm d}x .
\ee

\section{Proof of the Theorem}

\subsection{}

First of all (see \cite{7}, p. 79)
\bdis
Z(t)=e^{i\theta(t)}\zeta\left(\frac{1}{2}+it\right) \ \Rightarrow\ |Z(t)|=\left|\zeta\left(\frac{1}{2}+it\right) \right| .
\edis
Putting $\Delta=4$ into (\ref{f2.5}) we obtain
\be \label{f3.1}
\int_T^{T+U}Z^4\left[\frac{\varphi(t)}{2}\right]\hat{Z}^2(t){\rm d}t=2\int_{\frac{1}{2}\varphi(T)}^{\frac{1}{2}\varphi(T+U)}Z^4(x){\rm d}x,
\ee
i.e. we have to consider the integral
\be \label{f3.2}
\int_{\frac{1}{2}\varphi(T)}^{\frac{1}{2}\varphi(T+U)}Z^4(t){\rm d}t=\int_0^{\frac{1}{2}\varphi(T+U)}Z^4(t){\rm d}t-\int_0^{\frac{1}{2}}
Z^4(t){\rm d}t .
\ee

\subsection{}

Let us remind the Ingham - Heath-Brown formula (see \cite{2}, p. 129)
\be \label{f3.3}
\int_0^T Z^4(t){\rm d}t=T\sum_{k=0}^4 C_k(\ln T)^{4-k}+\mcal{O}(T^{7/8+\epsilon}) ,
\ee
which improved the Ingham formula (see \cite{1}, p. 277, \cite{7}, p. 125)
\bdis
\int_0^TZ^4(t){\rm d}t=\frac{1}{2\pi^2}T\ln^4 T+\mcal{O}(T\ln^3 T);\ C_0=\frac{1}{2\pi^2} .
\edis
Hence, by (\ref{f3.2}), (\ref{f3.3}) we have
\begin{eqnarray} \label{f3.4}
& &
\int_{\frac{1}{2}\varphi(T)}^{\frac{1}{2}\varphi(T+U)}Z^4(t){\rm d}t= \\
& &
=\left\{ \frac{\varphi(t)}{2}\sum_{k=0}^4 C_k\left(\frac{\varphi(t)}{2}\right)^{4-k}\right\}_{t=T}^{t=T+U}+\mcal{O}(T^{7/8+\epsilon})= \nonumber \\
& &
=\sum_{k=0}^4C_kV_k+\mcal{O}(T^{7/8+\epsilon}),\quad U=T^{7/8+2\epsilon} , \nonumber
\end{eqnarray}
where
\be \label{f3.5}
V_k=\frac{1}{2}\varphi(T+U)\ln^{4-k}\frac{\varphi(T+U)}{2}-\frac{1}{2}\varphi(T)\ln^{4-k}\frac{\varphi(T)}{2} .
\ee

\subsection{}

By the Ingham formula (see \cite{1}, p. 294, \cite{7}, p. 120)
\be \label{f3.6}
\int_0^T Z^2(t){\rm d}t=T\ln T+(2c-1-\ln 2\pi)T+\mcal{O}(T^{1/2}\ln T)
\ee
we have in the case $U=T^{7/8+2\epsilon}$
\be \label{f3.7}
\int_T^{T+U}Z^2(t){\rm d}t=U\ln T+(2c-\ln 2\pi)U+\mcal{O}(T^{7/8+\epsilon}) .
\ee
Next we have (see \cite{5}, (1.2))
\be \label{f3.8}
\int_T^{T+U}Z^2(t){\rm d}t=U\ln (T)\tan[\alpha(T,U)]\left\{ 1+\mcal{O}\left(\frac{\ln\ln T}{\ln T}\right)\right\} .
\ee
Comparing formulae (\ref{f3.7}) and (\ref{f3.8}) we obtain
\be \label{f3.9}
\frac{\varphi(T+U)-\varphi(T)}{2U}=\tan[\alpha(T,U)]=1+\mcal{O}\left(\frac{\ln \ln T}{\ln T}\right) .
\ee

\subsection{}

Since (see (\ref{f1.2}))
\be \label{f3.10}
T-\frac{1}{2}\varphi(T)\sim (1-c)\pi(T) \ \Rightarrow\ T\sim \frac{1}{2}\varphi(T) ,
\ee
it follows that
\begin{eqnarray} \label{f3.11}
& &
\ln\frac{\varphi(T+U)}{2}-\ln\frac{\varphi(T)}{2}=\ln\frac{\varphi(T+U)}{\varphi(T)}=\\
& &
=\ln\left[ 1+U\frac{\varphi(T+U)-\varphi(T)}{2U}\frac{1}{\varphi(T)}\right]=\ln\left[ 1+\mcal{O}\left(\frac{U}{T}\right)\right]= \nonumber \\
& &
=\mcal{O}\left(\frac{U}{T}\right)=\mcal{O}\left(\frac{1}{T^{1/8-2\epsilon}}\right) \nonumber ,
\end{eqnarray}
(see (\ref{f3.9}), (\ref{f3.10})). Hence, by (\ref{f3.11})
\begin{eqnarray} \label{f3.12}
& &
\ln\frac{\varphi(T+U)}{2}=\ln\frac{\varphi(T)}{2}+\mcal{O}\left(\frac{1}{T^{1/8-2\epsilon}}\right) , \\
& &
\ln^4\frac{\varphi(T+U)}{2}=\ln^4\frac{\varphi(T)}{2}+\mcal{O}\left(\frac{\ln^3 T}{T^{1/8-2\epsilon}}\right). \nonumber
\end{eqnarray}

\subsection{}

Next we have (see (\ref{f3.5}), (\ref{f3.9}) and (\ref{f3.12}))
\begin{eqnarray} \label{f3.13}
& &
V_0=U\frac{\varphi(T+U)-\varphi(T)}{2U}\ln^4\frac{\varphi(T+U)}{2}+ \\
& &
+\frac{1}{2}\varphi(T)\left[ \ln \frac{\varphi(T+U)}{2}-\ln\frac{\varphi(T)}{2}\right]\left[\ln^3\frac{\varphi(T+U)}{2}+\dots \right]= \nonumber \\
& &
U\left\{ 1+\mcal{O}\left( \frac{\ln\ln T}{\ln T}\right)\right\}\ln^4\frac{\varphi(T)}{2}
\left\{ 1+\mcal{O}\left(\frac{1}{T^{1/8-2\epsilon}}\right)\right\}+ \nonumber \\
& &
+\frac{\varphi(T)}{2}\mcal{O}\left(\frac{1}{T^{1/8-2\epsilon}}\right)\mcal{O}(\ln^3 T)= \nonumber \\
& &
=U\ln^4\frac{\varphi(T)}{2}+\mcal{O}(U\ln^3 T\ln\ln T)\sim U\ln^4 T . \nonumber
\end{eqnarray}
Similarly,
\be \label{f3.14}
V_l=\mcal{O}(U\ln^{4-l} T),\ l=1,2,3,4 \ .
\ee
Then from (\ref{f3.4}) by (\ref{f3.13}) and (\ref{f3.14}) the asymptotic formula
\be \label{f3.15}
\int_{\frac{1}{2}\varphi(T)}^{\frac{1}{2}\varphi(T+U)}Z^4(t){\rm d}t\sim C_0U\ln^4 T
\ee
follows, i.e. (see (\ref{f3.1}), (\ref{f3.15}))
\be \label{f3.16}
\int_T^{T+U}Z^4\left[\frac{\varphi(t)}{2}\right]\hat{Z}^2(t){\rm d}t\sim 2C_0U\ln^4 T .
\ee

\subsection{}

Next, using the mean-value theorem the left-hand side of (\ref{f3.16}) reduces to
\be\label{f3.17}
\frac{2}{\left\{ 1+\mcal{O}\left( \frac{\ln\ln \tau}{\ln \tau}\right)\right\}\ln\tau}\int_T^{T+U}Z^4\left[\frac{\varphi(t)}{2}\right]Z^2(t){\rm d}t,\
\tau\in (T,T+U) .
\ee
Since $\ln\tau\sim\ln T,\ T\to\infty$ then, from (3.16) by (3.17) the formula (1.1) follows.

\begin{remark}
The proof of (1.8) is similar to (3.1)-(3.17).
\end{remark}

\section{Concluding remarks}

Let

\bdis
S(t)=\frac{1}{\pi}\arg\zeta(1/2+it),\ S_1(T)=\int_0^T S(t){\rm d}t .
\edis
In one of my next papers the following correlation formulae
\begin{eqnarray*}
& &
\int_T^{T+U}\left\{S\left[\frac{\varphi(t)}{2}\right]\right\}^{2k}\left|\zeta(1/2+it)\right|^2{\rm d}t\sim
\frac{(2k)!}{k!(2\pi)^{2k}}U\ln T(\ln\ln T)^k, \\
& &
\int_T^{T+U}\left\{S_1\left[\frac{\varphi(t)}{2}\right]\right\}^{2k}\left|\zeta(1/2+it)\right|^2{\rm d}t\sim d_n U\ln T, \ \dots \ ,
\end{eqnarray*}
where $U=T^{1/2+\epsilon}$, will be studied.

\thanks{I would like to thank Michal Demetrian for helping me with the electronic version of this work.}

\end{document}